\documentclass[11pt]{amsart}
\usepackage{amsmath}

\newtheorem{theorem}{Theorem}[section]
\newtheorem{proposition}[theorem]{Proposition}
\newtheorem{lemma}[theorem]{Lemma}
\newtheorem{corollary}[theorem]{Corollary}

\theoremstyle{definition}

\theoremstyle{remark}
\newtheorem{remark}[theorem]{Remark}

\numberwithin{equation}{section}

\newcommand{\be}{\begin{equation}}
\newcommand{\ee}{\end{equation}}

\newcommand{\bbC}{{\mathbb C}}
\newcommand{\bbZ}{{\mathbb Z}}
\newcommand{\bbR}{{\mathbb R}}

\newcommand{\calY}{{\mathcal Y}}

\newcommand{\calT}{{\mathcal T}}
\newcommand{\calB}{{\mathcal B}}

\newcommand{\calU}{{\mathcal U}}
\newcommand{\calQ}{{\mathcal Q}}

\newcommand{\calL}{{\mathcal L}}

\newcommand{\calP}{{\mathcal P}}
\newcommand{\calH}{{\mathcal H}}
\newcommand{\calE}{{\mathcal E}}

\newcommand{\calR}{{\mathcal R}}
\newcommand{\calS}{{\mathcal S}}
\newcommand{\calW}{{\mathcal W}}

\newcommand{\calJ}{{\mathcal J}}

\newcommand{\fraku}{{\mathfrak u}}
\newcommand{\frakp}{{\mathfrak p}}

\DeclareMathOperator{\tr}{Tr}

\newcommand{\inner}[2]{\langle#1,#2\rangle}
\newcommand{\ket}[1]{|#1\rangle}
\newcommand{\bra}[1]{\langle #1 |}
\newcommand{\norm}[1]{\lVert#1\rVert}

\newcommand{\del}{\partial}

\newcommand{\Pinj}{\Pi_{N,J}}

\begin{document}

\title{Remarks on the naturality of quantization}
\author{T. Foth}
\address{Department of Mathematics\\
University of Western Ontario\\London, Ontario, Canada N6A 5B7}
\email{tfoth@uwo.ca}
\author{A. Uribe}
\address{Mathematics Department and Michigan Center for Theoretical Physics\\
University of Michigan\\Ann Arbor, Michigan 48109}
\email{uribe@umich.edu}
\thanks{A.U. supported in part by NSF grant DMS-0070690.}

\date{\today}
\begin{abstract}
Hamiltonian quantization of an integral compact symplectic
manifold $M$ depends on a choice of compatible almost-complex
structure, $J$.  For open sets $\calU$ in the set of compatible
almost complex structures and small enough values of Planck's constant,
the Hilbert spaces, $\calH_J$, of the quantization form a
bundle over $\calU$ with a natural ($L^2$) connection.
In this paper we examine the dependence of the
Hilbert spaces on the choice of $J$, by computing the
semi-classical limit of the curvature of this connection.
We also show that parallel transport provides a link between
the action of the group Symp$(M)$ of symplectic transformations
of $M$ and the
Schr\"odinger equation.
\end{abstract}

\maketitle
\tableofcontents

\newcommand{\symp}{\text{Symp}{^0}}

\section{Introduction}

Let $(M,\omega)$ be a compact integral symplectic manifold.
In this paper we consider the issue of the Hamiltonian quantization
of $(M,\omega)$, by which we mean the assignment of Hilbert
spaces to $M$ (depending on Planck's constant),
and self-adjoint
operators to real functions on $M$ inducing a deformation quantization
of $C^\infty(M)$. 
In recent years methods for quantizing $M$ in this sense
have been proposed
that involve the choice of:
\begin{enumerate}
\item A prequantum line bundle 
$L\to M$, i.e.\  a Hermitian line bundle
with connection $\nabla$ whose curvature is $\omega$, and
\item A compatible almost complex structure, $J$, on $M$
(instead of a polarization as in traditional geometric quantization).
\end{enumerate} 
Compatibility means that $J$ preserves $\omega$ and
that the symmetric bilinear form,
\[
g(\cdot , \cdot) = \omega(\cdot ,J(\cdot))
\]
is positive definite. Given a compatible $J$, one can define spaces 
\begin{equation}\label{uno}
\calE_{J,N}\subset L^2(M, L^{\otimes N})
\end{equation}
for all $N$ sufficiently large, see \cite{Ch}, \cite{Gu},\cite{BU}.  
(In case $(M,J,\omega)$ is a K\"ahler manifold,
$\calE_{J,N}$ is the space of  holomorphic sections of $L^{\otimes
N}$.)
Then, to a given classical Hamiltonian $H\in C^\infty(M)$, one
associates a sequence of operators $\text{Op}^N(H)$ on the
$\calE_{J,N}$ by the formula:
\[
\text{Op}^N(H) (\psi) = \Pinj \Bigl(-i\nabla_{\xi_H} + H \Bigr) (\psi),
\]
where $\nabla_{\xi_H}$ is covariant differentiation with respect
to the Hamilton vector field $\xi_H$ of $H$, and $\Pinj$ is the
orthogonal projection
\[
\Pinj: L^2(M, L^{\otimes N}) \to \calE_{J,N}.
\]
The association $H\mapsto \{\text{Op}^N(H)\}_N$ defines a star product
on $C^\infty (M)$.  More specifically, there exists a sequence of
bi-differential operators on $M$, $B_j$, such that
$\forall H,G\in C^\infty(M)$
\[
\text{Op}^N(H)\circ\text{Op}^N(G) \sim
\text{Op}^N(HG) + \sum_{j=1}^\infty N^{-j}\text{Op}^N(B_j(H,G)),
\]
and
\[
H\ast G := HG + \sum_{j=1}^\infty N^{-j}B_j(H,G)
\]
is a star product where $N=1/\hbar$ is the inverse of Planck's constant.
(The asymptotics above are in any operator norm.)

As we mentioned,
in case $(M,J,\omega)$ is a K\"ahler manifold, 
$\calE_{J,N}$ is the space of  holomorphic sections of $L^{\otimes
N}$.  If $(M,J,\omega)$ is not K\"ahler then one needs to use a notion
of ``nearly holomorphic" section.  We will not review here what this
notion means, but simply refer to {\it op.\ cit.\ }for details.
Suffice it to say that, in the non- integrable case, the
sequence of projections $\{ \Pinj\}_N$ has the same structure
as in the integrable case asymptotically as $N\to\infty$.  
(For the precise assumption see the statement of the main theorem below.)

It is natural to ask to what extent the quantization depends on the
choices of $(L,\nabla)$ and $J$.  The various choices of $(L,\nabla)$
are parametrized by the isomorphism classes of flat line bundles,
i.e.\ by $H^1(M,S^1)$.  In this paper we consider $(L,\nabla)$ and its
Hermitian structure fixed, and concentrate on the dependence on $J$.

Thus let $\calJ$ denote the infinite-dimensional space
of all compatible almost complex structures, and fix a
quantization scheme as in \cite{Gu}. One then gets
vector bundles
\[
\calE_N\to\calJ
\]
with fibers $\calE_{J,N}$ where $N$ is a large integer.  Actually,
the $\calE_{N,J}$ are defined only for large $N$ which can be chosen uniformly
with respect to $J$ in an open set, $\calU\subset\calJ$.  More precisely:  for
every $J\in\calJ$ there exists a neighborhood $\calU\ni J$ and an $N_0$ such that
for all $N>N_0$ $\calE_{J',N}$ is well-defined for $J'\in\calU$.
(This can be proved by a vanishing theorem using a Lichnerowicz' formula
and curvature estimates; the open set is in the $C^2$ topology to control
the curvature.) 
  
From the point of view of deformation quantization, all these quantizations are
isomorphic.  Specifically, for any $J$, $J'$ and $N$ large
enough there exist unitary operators 
$U_N:\calE_{J,N}\to \calE_{J',N}$ that induce isomorphisms of the 
deformation quantizations associated to $\{\calE_{J,N}\}_N$.
(See \S 4.2 in \cite{Ch}.)
However, one can ask whether the quantum
states themselves (i.e. rays in the Hilbert spaces $\calE_{J,N}$) are
in some sense independent of $J$.  This question was raised by Richard
Montgomery and Victor Ginzburg in \cite{GM}, in the following more
precise form: Is there a natural projectively flat connection on the
bundles $\calE_N$?  They point out that, since $L^2(M,L^{\otimes N})$ is
independent of $J$, (\ref{uno}) gives us a natural connection on each
of the bundles $\calE_N$, and ask whether this
connection is projectively flat, at least in the semi-classical limit
$N\to\infty$.  Specifically,
the natural connection on $\calE_N\to\calJ$ is defined as follows: 
Let $J(t)$ be a smooth curve on $\calJ$ and a
lift, $\psi_t \in\calE_{J(t),N}$ to the bundle $\calE_{N}$. 
Then the covariant
derivative of this lift is defined as:
\begin{equation}
D^N_{\dot{J}} \psi_t := \Pi_{N,J(t)}\Bigl(  \frac{d}{dt}\psi_t\Bigr).
\end{equation}
In this equation the derivative $\frac{d}{dt}\psi_t$ takes place in the $J$-independent
Hilbert space $L^2(M, L^{\otimes N})$, and $\Pinj: L^2(M, L^{\otimes N}) \to \calE_{J,N}$ is the orthogonal projection.

In this note we compute the large $N$ asymptotic behavior of the
curvature of $D^N$ and study some specific examples.  We will
see that the limit is never central.  
We note that Andr\'es Vi\~na
previously considered the issue of the projective flatness of this connection
(which he described somewhat differently), \cite{Vi}, and proved that
the curvature is non-zero for fixed $N$.  

It is known that for some special $M$ there
exist projectively flat submanifolds $\calY\subset\calJ$.  The prime examples
of this phenomenon are:
\begin{enumerate}
\item $M=\bbR^{2n}$ and $\calY$ is the space of {\em linear}
complex structures.  We review this example in the Appendix.
\item $M$ is a moduli space of flat connections on a closed surface, $\Sigma$,
and $\calY$ is the space of complex structures arising from complex structures
on $\Sigma$.
\end{enumerate}
The problem of finding other examples of $(M,\calY)$
appears to be quite challenging.

\bigskip\noindent
{\bf The semi-classical limit of the curvature.}
\medskip

In order to describe the limit of the curvature on $\calE_N\to\calJ$
as $N\to\infty$, we need to have a description of the tangent spaces,
$T_J\calJ$, for $J\in\calJ$.  
Let $J(t)$ be a
one-parameter family of almost complex structures on $M$.  
For each $x\in M$, from
the identity: $J_x^2 = -I$, one finds that the endomorphism
$A_x = \frac{d\ }{dt}J_x(t)|_{t=0}$ satisfies: $J_x(0)A_x+A_xJ_x(0)=0$.
In addition, if the $J(t)$ are compatible with $\omega$, one has that for
all $x$,
$\omega_x(A_x(\cdot),J_x(\cdot))+\omega_x(J_x(\cdot),A_x(\cdot))= 0$.  
Therefore
\begin{equation}\label{dos}
T_J\calJ =
\{\, A\in \text{End }(TM)\;;\;J\circ A+A\circ J=0\, 
\text{ and } \omega(A(\cdot),J(\cdot))+\omega(J(\cdot),A(\cdot))= 0 \}.
\end{equation}

Let $L\to M$ be as above.  Let $J\in\calJ$ and
$A,\; B\in T_J\calJ$.  Then the curvature of the natural
connection on $\calE_N$ evaluated on $(A,B)$,
$\Upsilon^{(N)}_J(A,B)$, is an operator
\[
\Upsilon^{(N)}_J(A,B): \calE_{J,N}\to\calE_{J,N}.
\]
The semi-classical limit of the curvature is described by the following:
\begin{theorem}\label{Main}
Assume that for each $J\in\calJ$ the sequence of orthogonal projectors,
$\Pi_{J,N}: L^2(M,L^{\otimes N})\to\calE_{J,N}$, is associated to a
Fourier integral operator of Hermite type in the unit circle bundle,
$Z\subset L^*$, symbolically identical to the Szeg\"o projector and
depending smoothly on $J$.
Then, given variations of complex structures $A,B$ as above, the
sequence of operators, $\Upsilon^{(N)}_J(A,B)$, is asymptotic to the
operator $T_\chi^{(N)} = \text{Op}^N(\chi)$ with symbol
\[
\chi_{A,B}(x) = \tr(A_xJ_xB_x).
\]
More precisely,
\[
\frac{1}{\text{dim }\calE_{J,N}}\,
\norm{\Upsilon^{(N)}_J(A,B) - T^{(N)}_\chi}_2^2 = O(1/N)
\]
as $N\to\infty$, where $\norm{A}_2^2= \tr{A^*A}$.
\end{theorem}

\medskip\noindent{\bf Remarks:}

\begin{enumerate}

\item The technical assumption on the $(\Pinj)$
holds if $M$ is K\"ahler and $\calJ$ is the
space of integrable compatible complex structures on $M$.
\item The function $\chi_{A,B}$ is general, in the following
sense: for any
$J$ and $\chi\in C^\infty(M)$, one can find $A$ and $B$ such that $\chi
= \chi_{A,B}$.
\item We will actually show that the sequence of operators,
$\{\Upsilon^{(N)}_J(A,B)\}$, constitutes a Fourier integral operator
of Hermite type, see Proposition \ref{Bcurv}.
\item A corollary of Theorem \ref{Main} is that {\em the curvature of
$\calE_N\to\calJ$ is highly non-central}.  Perhaps the simplest way to
see this is to observe that, by the Szeg\"o limit theorem for
Toeplitz operators and Theorem \ref{Main}, the normalized spectral
measure of $\Upsilon^{(N)}_J(A,B)$ converges weakly to the push-forward
of Liouville's measure on $M$ by the function $\chi_{A,B}$.
Centrality, in the semi-classical limit, would amount to the function
$\chi_{A,B}$ being {\em constant}.  On the other hand the curvature
preserves micro-support: As a Fourier integral operator its relation
can be identified with the diagonal on $M\times M$.
\end{enumerate}

\bigskip
As Donaldson points out in \cite{Do}, $\calJ$ is an infinite-dimensional
K\"ahler manifold.  With our description of $T\calJ$, the 
K\"ahler form, $\Omega$, on $\calJ$, is given by
\[
\Omega(A,B) = \int_M \tr(A_xJ_xB_x) \omega^n_x.
\]
On the other hand, a consequence of the previous Theorem is that
this is the limit of the trace of
$\frac{1}{\text{dim }\calE_{J,N}}\Upsilon^{(N)}_J(A,B)$.
Therefore, we obtain:
\begin{corollary}
The curvature of $\bigwedge^{\text{\tiny top}}\calE_N\to\calJ$
is asymptotic to a universal constant times $\text{vol}(M)N^n\Omega$,
where $n=\text{dim}(M)/2$.
\end{corollary}
Thus the metric connection on the determinant line bundle
$\bigwedge^{\text{\tiny top}}\calE_N\to\calJ$ is asymptotic
to a connection that quantizes $(\calJ, N\Omega)$. 

\bigskip\noindent
{\bf Parallel transport and the Schr\"odinger equation.}
\medskip

\newcommand{\vt}{v_t}

We now turn to the action of the group
$\symp$, of Hamiltonian diffeomorphisms of $M$,
on the bundles $\calE_N\to\calJ$ by the
``pull-back" operation.  We will then compare this action with
parallel transport along trajectories of $\symp$ on $\calJ$.
Our calculations will be independent of $N$, so we will suppress it
from the notations.

Let $H_t: M\to\bbR$ be a time-dependent Hamiltonian, and let
$\phi_t:M\to M$ be its associated Hamilton flow.  Thus $\phi_0$ is
the identity and if we let
$\xi_t$ be the vector field on $M$ given by:
\[
d(\phi_t)_x(\xi_t)_x= \frac{d}{dt}\phi_t(x),
\]
then $\iota_{\xi_t} \omega = -dH_t$.
The one-parameter family of symplectic transformations of $M$,
$\phi_t$, can be lifted to the bundle $L\to M$.  Specifically,
if $Z\subset L^*$ denotes the unit circle bundle, let us
define a time-dependent vector field on $Z$, $\eta_t$, by:
\[
\eta_t = \tilde \xi_t+ H_t\del_\theta,
\]
where $\tilde\xi$ is the horizontal lift of $\xi$ with respect to the
given connection and $\del_\theta$ is the generator of the $S^1$ action
on $Z$.  We let $f_t:Z\to Z$ the one-parameter family of
diffeomorphisms generated by $\eta_t$.  It is clear that for each $t$
$f_t$ commutes with the circle action, and therefore it defines an
automorphism of the bundle $L$ which we continue to denote by $f_t$.
We will denote by $V_t: L^2(M, L)\to L^2(M,L)$ the (unitary) operator
of pull-back with respect to $f_t$.  As is known, one has:
\begin{equation}
\dot V_t\circ V_t^{-1} = \nabla_{\xi_t} + iH_t.
\end{equation}

\medskip
Let us now pick a base point, $J_0\in\calJ$, and let
$\calE_t$ denote the fiber of $\calE\to\calJ$ over
$J_t = \psi_t^*J_0 = d(\psi_t^{-1})\circ J_0 \circ d(\psi_t)$,
and by $\Pi_t$ the orthogonal projector
onto $\calE_t$.  In the holomorphic (K\"ahler) case it
is clear that $V_t$ maps $\calE_0$ to $\calE_t$ and, $V_t$
being unitary, we therefore have:
\begin{equation}
\Pi_t = V_t^{-1}\circ\Pi_0\circ V_t.
\end{equation}
We will assume that this identity holds as well for the
chosen almost-holomorphic quantization scheme (naturality
of the scheme).

\begin{proposition}\label{Sch}
Denote by
$\calP_t: L^2(M,L)\to L^2(M,L)$ the operator defined by
\[
\calP_t|_{\calE_0}: \calE_0\to\calE_t \quad\text{is parallel
transport along } \{J_t\}
\]
and $\calP_t|_{\calE_0^{\bot}} = 0$.
Let $A: L^2(M,L)\to L^2(M,L)$ denote the Schr\"odinger operator
\begin{equation}
A_t= \Pi_0 \circ\left(-i \nabla_{\xi_t} + H_t\right)\circ\Pi_0,
\end{equation}
and let $S_t: L^2(M,L)\to L^2(M,L)$ denote the
solution to the problem: $\frac{d}{dt}S_t = iA_tS_t$, $S_0 =\Pi_0$.
Then
\begin{equation}
\calP_t = V_t^{-1}\circ S_t .
\end{equation}
\end{proposition}
Thus the fundamental solution of the Schr\"odinger equation arises as the 
operator that intertwines parallel transport and the operation of
pull-back by the classical Hamilton flow.  (Notice, by the way, that
if $J(t)\equiv J(0)$ then $V_t = S_t$, i.e.\ the fundamental solution
{\em is} the pull-back operator.)

\section{Proof of Theorem \ref{Main} and  Proposition \ref{Sch}}

\subsection{Preliminaries on the bundles $\calE_N\to\calJ$}\label{Sec2a}

We begin with a general expression for the curvature of the natural
connection of a bundle, $\calE\to \calJ$, such that 
$\forall J\in\calJ$ the fiber, $\calE_J$, is a closed subspace of a fixed
Hilbert space, $\calH$.  It will be useful to introduce the following
notation: For each $J\in\calJ$, we let $\Pi_J:\calH\to\calE_J$ be the
orthogonal projection, and if $v\in T_J\calJ$ we let $\delta_v\Pi$ denote
the operator on $\calH$
\[
\forall\psi\in\calH \qquad \delta_v\Pi (\psi) = 
\frac{d}{dt} \Pi_{\gamma(t)}(\psi),
\] 
where $\gamma(t)$ is a curve on $\calJ$ whose velocity at time $t=0$
is $v$.

\begin{lemma}\label{Curv}
If $v,w\in T_J\calJ$, the curvature of the natural
connection on $\calE$ at $(J,v,w)$ is the operator
\begin{equation}\label{3a}
\Upsilon_J(v,w) =
\Pi_J [\delta_w\Pi\,,\,\delta_v\Pi] : \calE_J\to\calE_J .
\end{equation}
\end{lemma} 
\begin{proof}
For completeness we sketch a proof.  Let $V,W$ be two commuting
vector
fields on $\calJ$, and let $\psi$ be a section of $\calE_N$. 
We wish to compute
$\Upsilon(V,W)(\psi) = \nabla_V(\nabla_W \psi) - \nabla_W(\nabla_V \psi)$.
If we denote by $\del_V\psi(J)\in\calH_{N,J}$ the derivative of
$\psi$ with respect to $V(J)$ inside the ambient Hilbert space
$L^2(M,L^{\otimes N})$, then
\[
\Upsilon(V,W)(\psi) = \Pi\del_V(\Pi\del_W \psi)-\Pi\del_W(\Pi\del_V \psi).
\]
Leibniz rule gives the desired result.
\end{proof}

\bigskip

We now prove Proposition \ref{Sch}.  We claim that
the following equations, together with the initial condition
$\calP_0=\Pi_0$,  characterize $\{\calP_t\}$:
\begin{equation}\label{pt1}
\Pi_t\circ \dot{\calP_t} = 0\,,\quad
\dot\Pi_t\circ\calP_t=\dot\calP_t.
\end{equation}
The left-hand side of the first equation is the definition of the
$L^2$ covariant derivative.
The second follows by differentiating the identities
$\calP_t\circ\Pi_0 = \Pi_t\circ\calP_t$ and
$\calP_t\circ (I-\Pi_0) = 0$, which themselves follow from the
definition of $\calP$ outside $\calE_0$.

We will show that $\calQ_t:=V_t^{-1}\circ S_t$ satisfies the
equations (\ref{pt1}),
and since $\calQ_0 =\calP_0=\Pi_0$ then necessarily
$\calQ_t=\calP_t$.  This is a computation, using the
fact that
\[
\Pi_0\circ S_t = S_t = S_t\circ\Pi_0.
\]
From this it follows that
\begin{equation}
\dot\calQ_t = V_t^{-1}\Bigl( -I + \Pi_0\Bigr)\dot V_t V_t^{-1} S_t
\end{equation}
and therefore $\Pi_t\circ\dot\calQ_t = 0$.  To check the second
equation (\ref{pt1}), simply use that
\begin{equation}
\dot\Pi_t =-V_t^{-1}\dot V_t V_t^{-1} \Pi_0V_t +
 V_t^{-1}\Pi_0 \dot V_t ,
\end{equation}
and therefore
\begin{equation}
\dot\Pi_t\calQ_t = V_t^{-1}
\Bigl( -\dot V_t V_t^{-1}\Pi_0 S_t + \Pi_0 \dot V_t V_t^{-1} S_t\Bigr) =
\end{equation}
\begin{equation}
= V_t^{-1}\Bigl( -I + \Pi_0\Bigr) \dot V_t V_t^{-1} S_t
= \dot\calQ_t.
\end{equation}

\bigskip
We now turn to the proof of Theorem \ref{Main}.

\subsection{Preliminaries on the linear case}\label{SecLin}

Let $(V,\omega)$ be a symplectic vector space, and let 
 $\calL$ denote the set of linear compatible complex structures
on $V$.  It will be useful to have several different descriptions of $\calL$:

\newcommand{\SpV}{\text{Sp}(V)}
\begin{enumerate}
\item $\calL$ is the symmetric space for the group $G =\SpV$ of
linear symplectic transformations of $V$.
\item $\calL$ is a Grassmannian:
We can identify $\calL$ with the set of positive complex
lagrangian subspaces of the complexification $V\otimes\bbC$, by
mapping a complex structure to the subspace of vectors of type $(0,1)$
with respect to that structure.  
\end{enumerate}
We now investigate the consequences of these two points of view.

Let us pick  a compatible complex structure on $V$, $J_0$, which will serve as a base point:
If we let $\text{U}\subset\SpV$ be the isotropy subgroup of $J_0$ (the unitary group
of $(V,\omega, J_0)$) then we can identify:
\[
\calL \cong \SpV/\text{U}.
\]
At the Lie algebra level the choice of $J_0$ induces a Cartan decomposition 
(see the Appendix for a matrix description):
\[
\text{sp} = \fraku + \frakp
\]
Elements of $\text{sp}$ are the linear transformations $A:V\to V$ such that
\[
\omega(A(\cdot),\cdot) + \omega(\cdot, A(\cdot)) = 0.
\]
The elements of $\frakp$ (resp.\ $\fraku$) are those $A$ that, in addition,
satisfy:
\[
A\circ J_0 = -A\qquad \text{resp. } A\circ J_0 = A.
\]
For any $A\in\text{sp}$
the vector field $V\ni v\mapsto A(v)$ is associated to the Hamiltonian
\begin{equation}
q_A(v) := \frac{1}{2}\omega(A(v),v).
\end{equation}
Therefore, elements of $\text{sp}$ can also
be thought of as (arbitrary) quadratic forms on $V$ (the Hamiltonians).  From this
point of view $\frakp$ is the space of quadratic forms that are the real part of
$J_0$-complex quadratic forms.
In a similar spirit, $\fraku$ consists of $J_0$-Hermitian quadratic forms.  

The vector space $\frakp$ has a symplectic structure:  
$\Omega (A,B) = \tr (AJ_0B)$, where we now think again of $A,B\in\frakp$ as linear
transformations of $V$.  This structure extends to the ($\SpV$-invariant) symplectic
form of $\calL$.

Let $\calB$ be a quantization of $(V,\omega)$: to be specific, let us pick
\[
\calB =
\{\, f:V\to\bbC\;;\; f\ J_0\text{-holomorphic}\,\}
\cap L^2(V, \exp(-\omega(v,J_0(v))\,\omega^n),
\]
the Bargmann space.  This Hilbert space carries a (Heisenberg) representation
of the Heisenberg group of $V$.   A lagrangian subspace, $\Lambda\subset V\otimes\bbC$
can be thought of as a maximal abelian subalgebra of the Lie algebra
of the complex Heisenberg group of $V$, heis$(V)_\bbC$. 
If $\Lambda$ is positive, the subspace,
$\calR_{\Lambda} \subset\calB$,
which is the joint kernel of all the operators induced 
by vectors of a positive $\Lambda$ is one-dimensional.  The spaces $\calR_\Lambda$
are the fibers of a line bundle, $\calR\to\calL$.  Since all the fibers are subspaces of the
fixed Hilbert space, $\calB$, we can endow this bundle with the $L^2$ connection.
Below we will need:

\begin{lemma}\label{CurvL}
The curvature of $\calR\to\calL$ is the symplectic form of $\calL$.
\end{lemma}
\begin{proof} 
We sketch a proof for convenience, invoking the existence of the metaplectic
representation of $V$.  Thus let Mp$(V)$ be the metaplectic group of $V$ (the
double cover of $\SpV$ on which the determinant has a smooth square root), 
and let $\rho$ denote the representation of Mp$(V)$ on $\calB$.  By its definition,
$\rho$ induces an action on the bundle $\calR\to\calL$, where $\rho$ acts on
$\calL$ via $\SpV$.  More specifically, if for each $\Lambda\in\calL$
$P_\Lambda:\calB\to\calB$ denotes the orthogonal projection onto $\calR_\Lambda$,
then
\begin{equation}
\forall g\in \text{Mp}(V)\qquad \rho(g)\circ P_\Lambda\circ\rho(g)^{-1} = P_{g\cdot\Lambda}.
\end{equation}
Let $A,B\in\frakp$.  By the formula above, the variation of $P_\Lambda$ by the infinitesimal
action of $A$ is
\[
\delta_AP_\Lambda = [d\rho(A)\,,\,P_\Lambda].
\]
By Lemma \ref{Curv} and the formula above, the curvature at $J_0$ evaluated at
$A,B\in\frakp$ is
\begin{equation}\label{lcurv}
\tr P_{\Lambda}\, [ [d\rho(B)\,,P_\Lambda]\; , \; [d\rho(A)\, P_\Lambda]].
\end{equation}
where $\Lambda$ is the space of $J_0$-$(0,1)$ vectors.
(Taking the trace is here used to identify the curvature as an operator on a 1-dimensional
vector space with a scalar.)  A short calculation shows that (\ref{lcurv}) is equal to
\[
\inner{d\rho([A,B]) e_\Lambda}{e_\Lambda},
\]
where $e_\Lambda\in\calB$ is a unit vector in $\calR_\Lambda$. 
In the Appendix we show that this inner product is indeed the
symplectic form of $\calL$ evaluated at $(A,B)$.
\end{proof}

\subsection{Hermite FIOs}

We need now to be more explicit about the quantization scheme we are
assuming.  Let $Z\subset L^*$ be the unit circle bundle and $\alpha$
the connection form on $Z$.  Then $\alpha$ is a contact form on $Z$ and
$d\alpha^n\wedge\alpha$ is a volume form on $Z$. If we let
\[
L^2(Z) = \oplus_{k\in\bbZ}L^2(Z)_k
\]
be the decomposition of $L^2(Z)$ into $S^1$ isotypes, then one has a
tautological isomorphism:
\[
L^2(Z)_k\cong L^2(M,L^{\otimes k}).
\]
For each $J\in\calJ$, let
\[
\calH_J := \bigoplus_{N=1}^\infty\calE_{J,N}
\]
(Hilbert space direct sum), and let $\Pi_J: L^2(Z)\to\calH_J$
be the orthogonal projector.
Let us also introduce the manifolds
\[
\Sigma = \{\, (p, r\alpha_p)\;;\;p\in Z, r>0 \,\},
\]
which is a symplectic submanifold of $T^*Z$, and
\[
\Sigma\stackrel{\Delta}\times\Sigma = \{\,((p,\sigma)\,,\,(p,-\sigma))\;;\;
(p,\sigma)\in\Sigma\,\}
\]
which is an isotropic submanifold of $T^*Z\times T^*Z$.  
The assumptions on our quantization scheme are:

\begin{enumerate}
\item For each $J\in\calJ$, the orthogonal
projector $\Pi_J: L^2(Z)\to\calH_J$ is
in the class $I^{1/2}(Z\times Z,\Sigma\stackrel{\Delta}\times\Sigma)$
of Boutet de Monvel and Guillemin, \cite{BG}, with a symbol to be
described below.
\item The map $J\mapsto \Pi_J$ is differentiable. 
\end{enumerate}
These assumptions hold, for example, if $\calJ$ is the space of
integrable structures and $\calH_{J,N}$ is the space of holomorphic
sections of $L^{\otimes N}$.

We now need to recall the notion of symbol of $\Pi_J$.  Since $\Sigma$ is
a symplectic submanifold of $T^*Z$, the bundle $\calS\to\Sigma$ where
\[
\calS_\sigma = (T_\sigma\Sigma)^o/T_\sigma\Sigma
\]
(where ${}^o$ denotes the symplectic orthogonal) is a bundle of symplectic
vector spaces over $\Sigma$.  $\calS$ is called the symplectic normal bundle,
and there is a natural isomorphism (induced by the projection,
$\pi: Z\to M$, lifted to $TT^*Z$):
\[
\calS_\sigma = T_{\pi(\sigma)} M^{-}.
\]
We will denote by $\calW_\sigma$ a Hilbert space which carries the 
metaplectic representation of the metaplectic group of $\calS_\sigma$.
Then, the symbol of a Fourier integral operator in a class
$I^{\cdot}(Z\times Z, \Sigma\stackrel{\Delta}\times\Sigma)$
is an object of the following type:  For each $\sigma\in\Sigma$, the
symbol at $\sigma$ is an operator
\[
\calW_\sigma\to \calW_\sigma,
\]
(if we take into consideration the metaplectic structure
then we have to work with the half forms; however the half forms will play
no role in our calculations). 
In particular, the symbol of $\Pi_J$
is a rank-one projector onto a ground state, $e_{J,\sigma}\in\calW_\sigma$:
\[
\ket{e_{J,\sigma}}\bra{e_{J,\sigma}}.
\]

\subsection{The curvature as an operator on $Z$.}
Let $J\in\calJ$ and $A,B\in T_J\calJ$.  As mentioned,
for each $N$, the curvature of our connection on $\calE_N\to\calJ$,
evaluated at $(A,B)$ is an operator
$\Upsilon^{(N)}_J(A,B): \calE_{J,N}\to\calE_{J,N}$.
Let us define:
\[
\Upsilon_J(A,B):=
\bigoplus_{N=1}^\infty\Upsilon^{(N)}_J(A,B): \calH_{J}\to\calH_{J},
\]
and extend $\Upsilon_J(A,B)$ to be an operator on $L^2(Z)$ by
defining it to be zero on $\calH_J^\bot$. (Recall that
$\calH_J = \oplus_N\calE_{J,N}$.)  We will prove:

\begin{proposition}\label{Bcurv}
The operator $\Upsilon_J(A,B)$ is a Fourier integral operator of Hermite
type, in the class $I^{1/2}(Z\times Z, \Sigma\stackrel{\Delta}\times\Sigma)$,
and  with symbol equal to the function $\chi_{A,B}$, pulled back to $\Sigma$,
times the symbol of $\Pi_J$.
\end{proposition}
\begin{proof}
We begin by proving that $\Upsilon_J(A,B)$ is an FIO in the stated class.
Let $J(s)$  be a smooth curve on $\calJ$ such that
$J(s)|_{s=0} = J$ and $\frac{d\ }{ds}J(s)|_{s=0} = A$.  
Introduce the notation:
\[
\delta_A\Pi_J := \frac{d\ }{dt}\Pi_{J(t)}\,|_{t=0},
\]
and analogously for $\delta_B\Pi_J$.  Then, by Lemma \ref{Curv},
\[
\Upsilon_J(A,B) = \Pi_J\,[\delta_A\Pi_J\,,\,\delta_B\Pi_J]\Pi_J.
\]
By the smooth dependence of  $\Pi_J$ on $J$, $\delta_A\Pi_J$ and
$\delta_B\Pi_J$ are both in the class
$I^{1/2}(Z\times Z, \Sigma\stackrel{\Delta}\times\Sigma)$.
It follows directly from Theorem 9.6 of \cite{BG} that $\Upsilon_J(A,B)$
is in this class as well.

\medskip

It remains to compute the symbol of $\Upsilon_J(A,B)$.  Let
$\sigma\in\Sigma$.  Denote by $x\in M$ the projection of
$\sigma$, and let $V = T_xM$.  This is a symplectic vector
space, with symplectic form $\omega_x$ and a compatible
linear complex structure, $J_x$.  The values $A_x$ and $B_x$
of $A$ and $B$ at $x$ are infinitesimal variations of the
linear structure $J_x$. 

Let $\calB$ be a Hilbert space quantizing $V$; for example 
we can take $\calB$ to be the Bargmann space
\[
\calB_{(V,J_x,\omega_x)} =
\{\, f:V\to\bbC\;;\; f\ J_x\text{-holomorphic}\,\}
\cap L^2(V, \exp(-\omega_x(v,J_x(v))\,\omega_x^n),
\]
and then the symbol of $\Pi_J$ at $\sigma$ is projection,
$P:\calB\to\calB$, onto the line in $\calB$ spanned by 
$\ket{e_J,\sigma}$ a (normalized) constant function.
Let us denote by $\delta_AP$, resp.\ $\delta_BP$ the 
derivative of $P$ with respect to $A_x$, resp.\ $B_x$.
By the general symbol calculus of Hermite FIOs,
\cite{BG}, the calculation of the symbol of $\Upsilon_J(A,B)$ at 
$\sigma$ amounts to the calculation of the operator
\begin{equation}
P\circ [\delta_BP\,,\,\delta_AP]\circ P.
\end{equation}
Thus the symbol calculus reduces the problem to the linear case
that we discussed in \S \ref{SecLin}. 
\end{proof}

\subsection{The end of the proof}

With Proposition \ref{Bcurv} at hand it is a simple matter to
finish the proof of the main theorem.  We keep the same notation
and assumptions as in the statement of Theorem \ref{Main}, except
that we will suppress the subindex $J$, for simplicity.
Let $\Xi$ denote the horizontal lift to $Z$ of the Hamilton vector
field of $\chi$, and consider the first-order operator on $Z$
\[
Q = -i\,\Xi + \chi.
\]
Then $\Pi Q\Pi $ is an Hermite FIO in the same class and with the
same symbol as $\Upsilon(A,B)$, and therefore the difference
\[
\calR := \Upsilon (A,B) - \Pi Q\Pi
\]
is in the class $I^{0}(Z\times Z, \Sigma\stackrel{\Delta}\times\Sigma)$.
Furthermore, it is clear that $\calR$ commutes with the circle action
on $Z$ and maps $\calH$ to itself.  Let $\calR^N: \calE_N\to\calE_N$
be the restriction of $\calR$ to $\calH_N\cong\calE_N$.
We wish to obtain the asymptotics of 
\[
\tau_N:= \text{trace of}\  (\calR_N)^*\circ\calR_N = 
\norm{\Upsilon^{(N)}(A,B) - T^{(N)}_\chi}_2^2.
\]
By Theorem 9.6 of \cite{BG}, the operator $\calR^*\circ\calR$
is in the class $I^{-1/2}(Z\times Z, \Sigma\stackrel{\Delta}\times\Sigma)$.
This means that the Schwartz kernel of $\calR^*\circ\calR$ is given
locally (up to a smooth error) by an oscillatory integral of
with the same phase as that of the Schwartz kernel of $\Pi$, but
with an amplitude of order one less in the fibre variables
than that of the Schwartz kernel of $\Pi$.  Therefore, the
trace $\tau_N$ has an asymptotic expansion of order one less
than the trace of the projectors onto the $\calE_N$, i.e.\ 
$\tau_N$ has an asymptotic expansion of order dim$\calE_N -1$.

\hfill{Q.E.D}

\def\a{{\frak{a}}}
\def\C{{\Bbb C}}
\def\d{{\partial}}
\def\g{{\frak{g}}}
\def\H{{\frak{H}}}
\def\k{{\frak{k}}}
\def\L{{\frak{L}}}
\def\M{{\frak{M}}}
\def\n{{\frak{n}}}
\def\p{{\frak{p}}}
\def\R{{\Bbb R}}
\def\T{{\cal{T}}}
\def\w{{\wedge}}

\newtheorem{Th}{Theorem}[section]
\newtheorem{assum}[Th]{Assumption}
\newtheorem{cor}[Th]{Corollary}

\section{Example: Teichm\"uller Space.}

In this section we compute explicitly the symbol
of the operator in Theorem 2.3 in the case
when $M$ is a compact surface of genus $g>1$
and the vector bundle $\calE _N$ is restricted to
the Teichm\"uller space $T_g$ embedded into $\calJ$ as a particular
slice.

Let $\Sigma$ be a compact smooth surface of genus $g>1$.
Choose and fix a hyperbolic metric $\sigma(z,\bar z)dzd\bar z$
on $\Sigma$. Here $\sigma (z,\bar z)=-\frac{2}{(z-\bar z)^2}$.
Denote by $Q(C)$ the space of holomorphic quadratic
differentials on the hyperbolic Riemann surface    
$C:=(\Sigma,\sigma(z,\bar z)dzd\bar z)$. 
Recall that the Teichm\"uller space is 
$\calT _g=Met_{-1}(\Sigma)/Diff_0(\Sigma)$,
and $Met_{-1}(\Sigma)\cong Met(\Sigma)/Conf(\Sigma)$
parametrizes all complex structures on $\Sigma$.
Consider an explicit
realization of the Teichm\"uller space $\calT _g$ as the smooth family
of Riemannian metrics on $\Sigma$ 
\begin{equation}
ds^2=f[\Phi dz^2+\rho Edzd\bar z+\bar\Phi d\bar z^2],
\label{fammet}
\end{equation}
where $w=w(z,\bar z)$ is uniquely determined by the condition
that the map $(\Sigma,\sigma(z,\bar z)dzd\bar z)\to 
(\Sigma,\rho (w,\bar w)dwd\bar w)$ which is the identity
as the map $\Sigma\to\Sigma$ is harmonic \cite{Wolf}, \cite{Jost},
$$
\rho (w,\bar w):=-\frac{2}{(w-\bar w)^2}, \ \ \
E:=\frac{\d w}{\d z}\frac{\d\bar w}{\d\bar z}+ 
\frac{\d w}{\d \bar z}\frac{\d\bar w}{\d z},
$$
$$
f:=\frac{\sigma(z,\bar z)}{\rho (w,\bar w)
(\frac{\d w}{\d z}\frac{\d\bar w}{\d\bar z}-
\frac{\d w}{\d \bar z}\frac{\d\bar w}{\d z})}, \ \ \
\Phi(z)=\rho (w,\bar w)\frac{\d w}{\d z}\frac{\d\bar w}{\d z},
$$
$\Phi(z)dz^2\in Q(C)$, and this is an explicit bijection 
between $\calT _g$ and $Q(C)$. Note that $E>0$.
All metrics (\ref{fammet}) have the same area form 
$\sigma(z,\bar z)dRe(z)\w dIm(z)$, this is achieved via introducing
the conformal factor $f$ which
is not present in the parametrization of $\calT_g$ given in
\cite{Wolf} (p.456) or \cite{Jost} (p.16), (\ref{fammet})
provides an embedding $\iota:\calT_g\hookrightarrow Met(\Sigma)$,
note that $\iota(\calT_g)\cap Met_{-1}(\Sigma)$
consists of one point (corresponding to $C$).

Now pick a point on the slice $\iota(\calT_g)$
represented by $\Phi_0dz^2\in Q(C)$.
The corresponding metric is 
$f_0[\Phi_0dz^2+\rho_0 E_0 dzd\bar z+\bar \Phi_0d\bar z^2],$
denote
$$
g_0=f_0\begin{pmatrix}\Phi_0+\bar \Phi_0+\rho_0E_0 & i(\Phi_0-\bar\Phi_0) \cr
i(\Phi_0-\bar\Phi_0) & -\Phi_0-\bar \Phi_0+\rho_0E_0.
\end{pmatrix}
$$
Choose two one-parameter families of metrics on $C$: 
$f_t[\Phi_tdz^2+\rho_t E_t dzd\bar z+\bar \Phi_td\bar z^2]$
and $f_s[\Phi_sdz^2+\rho_s E_s dzd\bar z+\bar \Phi_sd\bar z^2]$
with $\Phi_t|_{t=0}=\Phi_0$, $\Phi_s|_{s=0}=\Phi_0$.
We have:
$$
g_{J(t)}=f_t\begin{pmatrix}
\Phi_t+\bar \Phi_t+\rho_tE_t & i(\Phi_t-\bar\Phi_t) \cr
i(\Phi_t-\bar\Phi_t) & -\Phi_t-\bar \Phi_t+\rho_tE_t
\end{pmatrix},
$$
$$
g_{J(s)}=f_s\begin{pmatrix}
\Phi_s+\bar \Phi_s+\rho_sE_s & i(\Phi_s-\bar\Phi_s) \cr
i(\Phi_s-\bar\Phi_s) & -\Phi_s-\bar \Phi_s+\rho_sE_s
\end{pmatrix},
$$
$g_{J(t)}|_{t=0}=g_{J(s)}|_{s=0}=g_0$.
Denote $u_1=\frac{dg_{J(t)}}{dt}|_{t=0}$, 
$u_2=\frac{dg_{J(s)}}{ds}|_{s=0}$. So $u_1,u_2\in T_{g_0}\iota(\calT_g)$.

\begin{proposition}
The symbol of the Toeplitz operator in the Theorem 2.3 is 
\begin{equation}
\langle u_1,u_2\rangle =-\frac{8}{\sigma \rho_0f_0E_0} Im
(\frac{d(f_t\Phi_t)}{dt}|_{t=0}\frac{d(f_s\bar\Phi_s)}{ds}|_{s=0}).
\label{symb}
\end{equation}
\end{proposition}
\begin{remark}
If $\Phi_0\equiv 0$ (i.e. $g_0$ is the hyperbolic metric 
$\sigma(z,\bar z)dz d\bar z$ of $C$) then $f_0\equiv 1$, $E_0\equiv 1$,
$\rho_0\equiv \sigma$, and the pull-back to $T_{[C]}\calT_g$ of the natural
symplectic form on $T_{g_0}\iota(\calT_g)$ obtained via integrating 
(\ref{symb}) over $\Sigma$ with respect to the area form of $g_0$ 
is, up to a constant factor, the $(1,1)$-form associated to
the Weil-Petersson metric on $\calT_g$.
Indeed, the expression (\ref{symb}) becomes $-\frac{8}{\sigma^2}
Im(\varphi_1\bar\varphi_2)$, where $\varphi_1=\frac{d\Phi_t}{dt}|_{t=0}$, 
$\varphi_2=\frac{d\Phi_s}{ds}|_{s=0}$, $\varphi_1(z)dz^2,\varphi_2dz^2$
are holomorphic quadratic differentials on $C$ representing
two variations of the complex structure on $C$ (we use 
that $T_{[C]}(\calT _g)\simeq Q(C)$,
and also $T_0Q(C)\simeq Q(C)$ because $Q(C)$ is a finite-dimensional
complex vector space).
\end{remark}
\noindent {\bf Proof of the Proposition. }
First consider $\R^2$ with the symplectic form $\omega=c dx\w dy$,
where $c$ is a non-zero constant. Let us denote the matrices
of the symplectic form $\omega$, a complex structure $J$,
and the Riemannian metric $g_J$ defined by $g_J(u,v):=\omega(u,Jv)$
by the same letters. Denote $J_0=\begin{pmatrix} 0 & -1 \cr 1 & 0\end{pmatrix}$.
We have: $\omega=c J_0^{-1}$, $g_J=\omega J$, in particular
$g_{J_0}=c\begin{pmatrix} 1 & 0 \cr 0 & 1\end{pmatrix}$.
Consider a 1-parameter family $J(t)$ such that $J(0)=J_0$.
Denote $X:=\frac{dJ(t)}{dt}|_{t=0}\in \p$,
where $sp(1,\R)\simeq sl(2,\R)=\k + \p$.
The Riemannian metric obtained from the Killing form is,
up to a positive factor, $B(X_1,X_2)=tr(X_1X_2)$.
The invariant complex structure on $\p$ is $X\mapsto Z X Z^{-1}$,
where $Z:=\frac{1}{\sqrt{2}}\begin{pmatrix}1 & 1 \cr -1 & 1\end{pmatrix}$
\cite{Helgason} p.323. Therefore the symplectic form on $\p$
is given by $(X_1,X_2)=-B(X_1,Z X_2 Z^{-1})=-tr(X_1 Z X_2 Z^{-1})$.
A tangent vector to the space of linear complex structures
on $\R^2$ viewed as the set of metrics with the area form $\omega$
is $\frac{dg_{J(t)}}{dt}|_{t=0}=\frac{d(\omega J(t))}{dt}|_{t=0}=\omega X
=cJ_0^{-1}X$.

Now pick a point $p\in \Sigma$ and apply the above considerations 
to $T_p\Sigma\simeq (\R^2,g_0|_{T_p\Sigma})$. 
Everything is assumed to be restricted to $T_p\Sigma$.
We have: $c=\sigma$, 
$g_0=\omega hJ_0h^{-1}=\sigma J_0^{-1}hJ_0h^{-1}$,
where 
$$
h=\frac{1}{\sqrt{2\sigma(\sigma+f_0\rho_0E_0)}}\begin{pmatrix}
-\sigma+f_0(\Phi_0+\bar\Phi_0-\rho_0E_0) & if_0(\Phi_0-\bar\Phi_0) \cr
if_0(\Phi_0-\bar\Phi_0) & -\sigma-f_0(\Phi_0+\bar\Phi_0+\rho_0E_0)\end{pmatrix},
$$
$h\in Sp(1,\R)$,
$$
X_1=h^{-1}\omega^{-1}\frac{dg_{J(t)}}{dt}|_{t=0}h=
\frac{1}{\sigma}h^{-1}J_0\frac{dg_{J(t)}}{dt}|_{t=0}h,
$$
$$
X_2=h^{-1}\omega^{-1}\frac{dg_{J(s)}}{ds}|_{s=0}h=
\frac{1}{\sigma}h^{-1}J_0\frac{dg_{J(s)}}{ds}|_{s=0}h,
$$
and
$$
\langle \frac{dg_{J(t)}}{dt}|_{t=0},\frac{dg_{J(s)}}{ds}|_{s=0}\rangle =
-\frac{1}{\sigma^2}
tr (h^{-1}J_0\frac{dg_{J(t)}}{dt}|_{t=0}h Z
h^{-1}J_0\frac{dg_{J(s)}}{ds}|_{s=0}h Z^{-1}).
$$
A straightforward computation results in (\ref{symb}).
{\bf End of the proof.}

\begin{appendix}

\newtheorem{prop}[Th]{Proposition}
\def \proof{{\noindent{\it Proof.\ \ }}}
\newtheorem{lem}[Th]{Lemma}

\section{The Euclidean case}

Consider $V=\R^{2n}$ with the standard symplectic form
$\omega =\sum_{j=1}^{n}dx_j\w dy_j=\frac{i}{2}\sum_{j=1}^{n}dz_j\w
d\bar z_j$, $x,y\in\R^n$, $z=x+iy\in\C^n$.
Denote $\sigma=\begin{pmatrix}0_n & -1_n \cr 1_n & 0_n\end{pmatrix}$.
We have: $\omega (u,v)=u^t\sigma ^t v$.
Notations: $G=Sp(n,\R)=\{ X\in GL(2n,\R)|X^t\sigma X=\sigma\}$,
$\g =sp(n,\R)=\{ X\in Mat(2n,\R)|X^t\sigma +\sigma X=0\}=
\{ \begin{pmatrix} A & B \cr C & -A^t\end{pmatrix}|B=B^t, C=C^t\}$,
a maximal compact subgroup of $G$ $K=Sp(n,\R)\cap U(2n)$,
Cartan automorphism $\Theta:\g\rightarrow \g$, 
$X\mapsto -X^t$, the Cartan
decomposition is 
\begin{equation}
\g=\k +\p , 
\label{cartan}
\end{equation}
where
$$
\k = sp(n,\R)\cap u(2n)=\{ X\in Mat(2n,\R)|X^t\sigma +\sigma X=0, X+X^t=0\} =
$$
$$
\{ \begin{pmatrix} A & B \cr -B & A\end{pmatrix}|A=-A^t, B=B^t\} =
\{ X\in\g |\Theta(X)=X\},
$$
and
$$
\p=\{ X\in\g |\Theta(X)=-X\} =
\{ \begin{pmatrix} A & B \cr B & -A\end{pmatrix}|A=A^t, B=B^t\}.
$$
Also we have: $G=PK$, where $P=\exp \p$, and $\g_{\C}=\k_{\C}+\p_{\C}$.

The decomposition of an arbitrary element $X\in \g$ according
to (\ref{cartan}) is given by
$$
X=\frac{1}{2}(X-X^t)+\frac{1}{2}(X+X^t).
$$
The algebra $\g$ is isomorphic to the algebra of quadratic Hamiltonians
(with the Poisson bracket). For $X\in\g$ the corresponding
quadratic form is 
$H(v)=\frac{1}{2}\omega (v,Xv)=\frac{1}{2}v^t\sigma ^t Xv$.
The Poisson bracket is given by
$\{ H_1,H_2\}(v) =\frac{1}{2}\omega(v,[X_2,X_1]v)=\omega (X_1v,X_2v)$.
For our purposes we would like to decompose $H(v)$ according
to (\ref{cartan}). We obtain:
$$
H(v)=v^t(\frac{1}{4}\sigma ^t (X-X^t)+\frac{1}{4}\sigma ^t (X+X^t))v=
$$
\begin{equation}
\frac{1}{4}Re(iz^t(A'+iB')\bar z+iz^t(A''-iB'')z),
\label{decomp}
\end{equation}
where $X=\begin{pmatrix} A & B \cr C & -A^t\end{pmatrix}\in \g$, $B=B^t$, $C=C^t$,
$A'=A-A^t$, $B'=B-C^t$, $A''=A+A^t$, $B''=B+C^t$. We observe:
$A'+iB'\in u(n)$, and $A''-iB''$ is a symmetric complex matrix.

The Hamiltonian vector field ${\xi _H}$ 
is determined by ${\xi_H}\lrcorner \omega =dH$, 
in local coordinates
$$
{\xi_H}=\sum_{j=1}^n(\frac{\d H}{\d y_j}\frac{\d}{\d x_j}
-\frac{\d H}{\d x_j}\frac{\d}{\d y_j} )=
2i\sum_{j=1}^n(\frac{\d H}{\d z_j}\frac{\d}{\d \bar z_j}
-\frac{\d H}{\d \bar z_j}\frac{\d}{\d z_j} ).
$$
We shall denote by $\L_{\xi_H}$ the  derivation on $C^{\infty}(V)$
induced by $\xi_H$, the Lie bracket is
$[\L_{\xi_{H_1}}, \L_{\xi_{H_2}}]= \L_{\xi_{H_1}}\L_{\xi_{H_2}}-
\L_{\xi_{H_2}}\L_{\xi_{H_1}}$.

The space $\calJ$ of all linear complex structures on $V$ compatible with
the standard symplectic form $\omega$ is  
$$
\calJ =\{ J\in GL(2n,\R)|J^2=-1, J^t\sigma J=\sigma, 
v^t\sigma ^t J v>0 \ for \ v\ne 0\}
$$
and is naturally identified with the bounded symmetric domain of type III:
$\calJ\cong Sp(n,\R)/U(n)$, $G$ acts on $\calJ$ by conjugation. Any complex structure
$J\in \calJ$ is obtained from the standard complex structure $\sigma$
as $J=g\sigma g^{-1}$ for $g\in Sp(n,\R)$, in fact $g$ can always be 
chosen in $P$.

A function $\psi :V\rightarrow \C$ is called {\it $J$-holomorphic} 
if $d\psi \circ J=id\psi$. 

For a positive integer $N$ define the Bargmann space
$$
H _{\sigma}^{(N)}=\{ f(z)e^{-\frac{N}{2}|z|^2}| 
f:\C^n\rightarrow\C \ is\ entire, \ 
\int |f(z)|^2e^{-N|z|^2}d\mu (z)<\infty\}=
$$
$$ 
L^2_{holo}(\C^n; e^{-N|z|^2}d\mu(z)),
$$
where $d\mu (z)$ is the standard Lebesgue measure.
It is in natural bijection with the space of $\sigma$-holomorphic 
sections of the line bundle $L^{\otimes N}$, where 
$L\rightarrow \C^n$ is the trivial line bundle
with Hermitian structure defined by $(s,s)=e^{-\frac{|z|^2}{2}}$, 
where $s$ is the unit section. $L$ is a quantizing
line bundle on $\C^n$ in the sense that $c_1(L)=\frac{1}{2\pi}[\omega]$.
Denote by $P$ the unit circle bundle in $L^*$.
Denote by $\pi=\pi_\sigma ^{(N)}$ the orthogonal projector 
$L^2(P)\rightarrow H_{\sigma}^{(N)}$. We recall that
$$
\pi (f(z,\bar z)e^{-\frac{N}{2}|z|^2})=
ce^{-\frac{N}{2}|z|^2}\int f(w,\bar w)e^{Nz\bar w}e^{-N|w|^2}d\mu (w),
$$
where $c=c(n,N)$ is a constant \cite{Bort}.

The following lemma is proved using integration by parts,
the proof is a short straightforward computation.
\begin{lem}
For $f(z)e^{-\frac{N}{2}|z|^2}\in H_\sigma^{(N)}$ 
the following identities hold:
$$
\pi (\bar z_s f(z)e^{-\frac{N}{2}|z|^2})=
\frac{1}{N}\frac{\d f}{\d z_s}e^{-\frac{N}{2}|z|^2},
$$
$$
\pi (\bar z_s\bar z_r f(z)e^{-\frac{N}{2}|z|^2})=
\frac{1}{N^2}\frac{\d ^2 f}{\d z_s\d z_r}e^{-\frac{N}{2}|z|^2}.
$$
\end{lem}

We observe: if $g\in G$, $J=g\sigma g^{-1}$, and $\psi$
is $\sigma$-holomorphic, then 
$d(\psi\circ g^{-1})\circ J=d\psi\circ \sigma\circ g^{-1}=
id\psi \circ g^{-1}=id(\psi\circ g^{-1})$,
so $\psi\circ g^{-1}$ is $J$-holomorphic. 
The map $\psi\mapsto \psi\circ g^{-1}$ establishes a unitary isomorphism
between $H_{\sigma}^{(N)}$ and the space $H_{J}^{(N)}$
of $J$-holomorphic sections of $L^{\otimes N}$. 
We shall drop $N$ from notation and write simply $\pi_J$.

Denote by $J_0$ be the fixed point of $K$ in $\calJ$.
\begin{prop}
Let $J_1(t)=g_1(t)\sigma g_1(t)^{-1}$ and 
$J_2(t)=g_2(t)\sigma g_2(t)^{-1}$ be two paths
in $\calJ$ such that $J_1(0)=J_2(0)=J_0$, $g_1(t)\in P$, $g_2(t)\in P$,
$g_1(0)=g_2(0)=I$.
Denote by $\tau_{g_1(t)}$ and $\tau_{g_2(t)}$ the corresponding
endomorphisms $V\rightarrow V$, also denote 
$\dot{\tau}_1=\dot{\tau}_{g_1(0)}$, 
$\dot{\tau}_2=\dot{\tau}_{g_2(0)}$.
The curvature of the connection at $J_0$ is
equal to 
$$
\pi [\dot{\tau}_2,\dot{\tau}_1]\pi -
[\pi\dot{\tau}_2\pi,\pi \dot{\tau}_1 \pi ].
$$
\end{prop}
\proof 
For $\pi_{J_1(t)}=\tau_{g_1(t)}\circ\pi \circ\tau_{g_1(t)}^{-1}$
we obtain:
$$
\dot{\pi}_{J_1(t)}=\dot{\tau}_{g_1(t)}\circ\pi\circ\tau_{g_1(t)}^{-1}-
\tau_{g_1(t)}\circ\pi\circ\tau_{g_1(t)}^{-1}\circ
\dot{\tau}_{g_1(t)}\circ\tau_{g_1(t)}^{-1},
$$
where $\dot{\tau}_{g(t)}(\varphi)(v)=\frac{d}{dt}\varphi(g(t)v)$,
and the hamiltonian operator $\dot{\tau_1}$
acts by $\dot{\tau_1}(\varphi)(v)=\frac{d}{dt}|_{t=0}\varphi(g(t)v)$.
So $\dot{\pi}_{J_1(t)}|_{t=0}=[\dot{\tau_1},\pi ]$,
and similarly $\dot{\pi}_{J_2(t)}|_{t=0}=[\dot{\tau_2},\pi ]$.
Finally
$$
Curv(\nabla)=\pi [[\dot{\tau_1},\pi ],
[\dot{\tau_2},\pi ]]\pi =
\pi (-\dot{\tau_1}\dot{\tau_2}+\dot{\tau_1}\pi \dot{\tau_2}+
\dot{\tau_2}\dot{\tau_1}-\dot{\tau_2}\pi\dot{\tau_1})\pi =
$$
$$
\pi [\dot{\tau}_2,\dot{\tau}_1]\pi -
[\pi\dot{\tau}_2\pi,\pi \dot{\tau}_1 \pi ].
$$
Q.E.D.

Denote $Curv_{H_1,H_2}=\pi [\L_{\xi_{H_2}},\L_{\xi_{H_1}}]\pi -
[\pi \L_{\xi_{H_2}}\pi, \pi \L_{\xi_{H_1}}\pi]$.
Decomposition (\ref{decomp}) shows that 
computation of curvature for $\dot{\tau}_1\in \p$, 
$\dot{\tau}_2\in\p$ is reduced to the following proposition.
\begin{prop}

\noindent (i) If $H_1=z_mz_l$, $H_2=z_rz_s$, then 
$Curv_{H_1,H_2}=0$,

\noindent  (ii) If $H_1=\bar z_m\bar z_l$, $H_2=\bar z_r\bar z_s$, then 
$Curv_{H_1,H_2}=0$,

\noindent (iii) If $H_1=z_mz_l$, $H_2=\bar z_r\bar z_s$, then 
$Curv_{H_1,H_2}=4(\delta_{mr}\delta_{ls}+\delta_{ms}\delta_{lr})$.
\label{curvvalue}
\end{prop}
\begin{remark}
These quadratic Hamiltonians are associated to elements 
of $\p_{\C}$.
\end{remark}
\proof (i) We have: 
$$
\xi_{H_1}=2i(z_l\frac{\d}{\d \bar z_m}+z_m\frac{\d}{\d \bar z_l}),
\xi_{H_2}=2i(z_r\frac{\d}{\d \bar z_s}+z_s\frac{\d}{\d \bar z_r}),
[\L_{\xi_{H_2}},\L_{\xi_{H_1}}]=0,
$$
$$
\pi\L_{\xi_{H_1}}\pi\L_{\xi_{H_2}}\pi (f(z)e^{-\frac{k}{2}|z|^2})=
-2ik\pi \L_{\xi_{H_1}}\pi (f(z)z_rz_se^{-\frac{k}{2}|z|^2})=
$$
$$
-4k^2f(z)z_rz_sz_lz_me^{-\frac{k}{2}|z|^2}.
$$
Similarly
$$
\pi\L_{\xi_{H_2}}\pi\L_{\xi_{H_1}}\pi (f(z)e^{-\frac{k}{2}|z|^2})=
-4k^2f(z)z_rz_sz_lz_me^{-\frac{k}{2}|z|^2},
$$
therefore $Curv_{H_1,H_2}=0$.

(ii)
$$
\xi_{H_1}=-2i(\bar z_l\frac{\d}{\d z_m}+\bar z_m\frac{\d}{\d z_l}),
\xi_{H_2}=-2i(\bar z_r\frac{\d}{\d z_s}+\bar z_s\frac{\d}{\d z_r}),
[\L_{\xi_{H_2}},\L_{\xi_{H_1}}]=0,
$$
$$
\pi\L_{\xi_{H_1}}\pi\L_{\xi_{H_2}}\pi (f(z)e^{-\frac{k}{2}|z|^2})=
-2i\pi \L_{\xi_{H_1}}\pi ((\bar z_s\frac{\d f}{\d z_r}+
\bar z_r\frac{\d f}{\d z_s}
-k\bar z_s\bar z_rf(z))e^{-\frac{k}{2}|z|^2})=
$$
$$
-\frac{2i}{k}\pi \L_{\xi_{H_1}}(\frac{\d ^2f}{\d z_s\d z_r}
e^{-\frac{k}{2}|z|^2})=
$$
$$
-\frac{4}{k}\pi ((\bar z_l \frac{\d}{\d z_m}\frac{\d ^2f}{\d z_s\d z_r}+
\bar z_m \frac{\d}{\d z_l}\frac{\d ^2f}{\d z_s\d z_r}
-k\bar z_l\bar z_m\frac{\d ^2f}{\d z_s\d z_r})e^{-\frac{k}{2}|z|^2})=
$$
$$
-\frac{4}{k^2}\frac{\d ^2}{\d z_l\d z_m}
\frac{\d ^2f}{\d z_s\d z_r} e^{-\frac{k}{2}|z|^2}.
$$
Computing $\pi\L_{\xi_2}\pi\L_{\xi_1}\pi (f(z)e^{-\frac{k}{2}|z|^2})$
we obtain the same expression, hence $Curv_{H_1,H_2}=0$.

(iii)
$$
\xi_{H_1}=2i(z_l\frac{\d}{\d \bar z_m}+z_m\frac{\d}{\d \bar z_l}),
\xi_{H_2}=-2i(\bar z_r\frac{\d}{\d z_s}+\bar z_s\frac{\d}{\d z_r}),
$$
$$
[\L_{\xi_{H_2}},\L_{\xi_{H_1}}]=4(
\delta_{sl}(-z_m\frac{\d}{\d z_r}+\bar z_r\frac{\d}{\d \bar z_m})+ 
\delta_{rl}(-z_m\frac{\d}{\d z_s}+\bar z_s\frac{\d}{\d \bar z_m}) )+
$$
$$
\delta_{sm}(-z_l\frac{\d}{\d z_r}+\bar z_r\frac{\d}{\d \bar z_l})+
\delta_{rm}(-z_l\frac{\d}{\d z_s}+\bar z_s\frac{\d}{\d \bar z_l})).
$$
We compute:
$$
\pi\L_{\xi_{H_1}}\pi\L_{\xi_{H_2}}\pi (f(z)e^{-\frac{k}{2}|z|^2})=
-2i\pi\L_{\xi_{H_1}}\pi ((\bar z_s\frac{\d f}{\d z_r} 
+\bar z_r\frac{\d f}{\d z_s}- 
kf(z)\bar z_s\bar z_r)e^{-\frac{k}{2}|z|^2})=
$$
$$
-\frac{2i}{k}\pi\L_{\xi_{H_1}}
(\frac{\d ^2f}{\d z_s\d z_r} e^{-\frac{k}{2}|z|^2})=
-4\pi (\frac{\d ^2f}{\d z_s\d z_r} z_lz_m
e^{-\frac{k}{2}|z|^2}),
$$
hence
\begin{equation}
\pi\L_{\xi_{H_1}}\pi\L_{\xi_{H_2}}\pi (f(z)e^{-\frac{k}{2}|z|^2})=
-4\frac{\d ^2f}{\d z_s\d z_r} z_lz_m
e^{-\frac{k}{2}|z|^2}.
\label{L1L2}
\end{equation}
Also
$$
\pi\L_{\xi_{H_2}}\pi\L_{\xi_{H_1}}\pi (f(z)e^{-\frac{k}{2}|z|^2})=
-2ik\pi\L_{\xi_{H_2}}\pi(f(z)z_lz_me^{-\frac{k}{2}|z|^2})=
$$
$$
-4k\pi ((\bar z_s\frac{\d (f(z)z_lz_m)}{\d z_r}+
\bar z_r\frac{\d (f(z)z_lz_m)}{\d z_s}-kf(z)z_lz_m\bar z_s\bar z_r)
e^{-\frac{k}{2}|z|^2}),
$$
so
\begin{equation}
\pi\L_{\xi_{H_2}}\pi\L_{\xi_{H_1}}\pi (f(z)e^{-\frac{k}{2}|z|^2})=
-4\frac{\d ^2 (f(z)z_lz_m)}{\d z_r\d z_s}e^{-\frac{k}{2}|z|^2}.
\label{L2L1}
\end{equation}
Also we have:
$$
\pi[\L_{\xi_{H_2}},\L_{\xi_{H_1}}]\pi (f(z)e^{-\frac{k}{2}|z|^2})=
$$
\begin{equation}
-4(z_l(\delta_{sm}\frac{\d f}{\d z_r}+\delta_{rm}\frac{\d f}{\d z_s})+
z_m(\delta_{ls}\frac{\d f}{\d z_r}+\delta_{rl}\frac{\d f}{\d z_s}))
e^{-\frac{k}{2}|z|^2}.
\label{commL1L2}
\end{equation}
Combining together (\ref{L1L2}), (\ref{L2L1}), and (\ref{commL1L2}),
we obtain:
$$
Curv_{H_1,H_2} (f(z)e^{-\frac{k}{2}|z|^2})=
4(\delta_{mr}\delta_{ls}+\delta_{ms}\delta_{lr})f(z)e^{-\frac{k}{2}|z|^2}.
$$
Q.E.D.

The subspace of the algebra of quadratic hamiltonians 
corresponding to $\p$ is generated (over $\R$) by 
$H^+_{m,l}:=z_mz_l+\bar z_m\bar z_l$
and  $H^-_{m,l}:=i(z_mz_l-\bar z_m\bar z_l)$, $1\le l\le n$,  $1\le m\le n$.
Denote by $\p^+$ (resp. $\p^-$) the subspace of $\p$ corresponding to the 
linear span of $H^+_{m,l}$ (resp. $H^-_{m,l}$) in the algebra of
quadratic hamiltonians. We have the decomposition $\p=\p^+\oplus \p^-$.
\begin{remark}
$\dim\g=2n^2+n$, $\dim \k=n^2$, $\dim\p =n^2+n$, 
$\dim \p^+=\dim\p^-=\frac{n(n+1)}{2}$, and $\dim \a= n$, where 
$\p=\a +\n$.
\end{remark}
We see that the curvature restricted to $\p^+$ (and to $\p^-$) is zero,
indeed, Proposition \ref{curvvalue} implies:
\begin{cor}
$$
Curv_{H^+_{m,l};H^-_{r,s}}=-8i(\delta_{mr}\delta_{ls}+\delta_{ms}\delta_{lr}),
$$
$$
Curv_{H^+_{m,l};H^+_{r,s}}=Curv_{H^-_{m,l};H^-_{r,s}}=0.
$$
\label{curv-explicit}
\end{cor}
\proof 
$$
Curv_{H^+_{m,l};H^-_{r,s}}=i(
Curv_{z_mz_l,z_rz_s}+Curv_{\bar z_m\bar z_l,z_rz_s}-
Curv_{z_mz_l,\bar z_r\bar z_s}-Curv_{\bar z_m\bar z_l,\bar z_r\bar z_s})=
$$
$$
i(-Curv_{z_rz_s,\bar z_m\bar z_l}-Curv_{z_mz_l,\bar z_r\bar z_s})=
-8i(\delta_{mr}\delta_{ls}+\delta_{ms}\delta_{lr})).
$$
Similarly
$$
Curv_{H^+_{m,l};H^+_{r,s}}=
Curv_{z_mz_l,z_rz_s}+Curv_{\bar z_m\bar z_l,z_rz_s}+
Curv_{z_mz_l,\bar z_r\bar z_s}+Curv_{\bar z_m\bar z_l,\bar z_r\bar z_s}=
$$
$$
-Curv_{z_rz_s,\bar z_m\bar z_l}+Curv_{z_mz_l,\bar z_r\bar z_s}=
-4(\delta_{mr}\delta_{ls}+\delta_{ms}\delta_{lr})
+4(\delta_{mr}\delta_{ls}+\delta_{ms}\delta_{lr})=0,
$$
and
$$
Curv_{H^-_{m,l};H^-_{r,s}}=-(
Curv_{z_mz_l,z_rz_s}-Curv_{\bar z_m\bar z_l,z_rz_s}-
Curv_{z_mz_l,\bar z_r\bar z_s}+Curv_{\bar z_m\bar z_l,\bar z_r\bar z_s})=
$$
$$
-(Curv_{z_rz_s,\bar z_m\bar z_l}-Curv_{z_mz_l,\bar z_r\bar z_s})=
-(4(\delta_{mr}\delta_{ls}+\delta_{ms}\delta_{lr})
-4(\delta_{mr}\delta_{ls}+\delta_{ms}\delta_{lr}))=0.
$$
Q.E.D.

\

\begin{cor}
If $J\in\calJ$, $u,v\in T_J\calJ$, and $q_1$, $q_2$ are the quadratic forms
representing $u,v$, then
the curvature of the connection (which is, once $u,v$ are fixed,
is an element of $End(H_J^{(k)})$) 
is equal to $i\Omega (q_1,q_2)Id$, where $Id$ denotes
the identity operator on the fiber $H_J^{(k)}$.   
\end{cor}
\proof
First we see that due to Corollary \ref{curv-explicit}
this is true if $J=J_0$ (the proof is straightforward verification).
The curvature and $\Omega$ are both
$G$-invariant, therefore the equality holds at every other point of $\calJ$.

Q.E.D.
\newpage
\section{Notation Index}

\addtolength{\baselineskip}{5pt}

\begin{itemize}
\item $(M,\omega)$: the given symplectic manifold.

\item $L\to M$: a pre-quantum line bundle of $M$.

\item $\nabla:$  The connection on $L\to M$.
\item $\calJ$: the infinite-dimensional manifold of 
almost complex structures.
\item $\Omega$: The K\"ahler form of $\calJ$.
\item $\calE_{J,N}\subset L^2(M,L^N)$: the Hilbert space associated to $J$.
\item $\calE_N\to\calJ$: the bundle of Hilbert spaces, with fibers the
$\calE_{N,J}$.
\item $D^N$:  The natural connection on $\calE_N\to\calJ$.
\item $\Upsilon^{(N)}_J(A,B): \calE_{J,N}\to\calE_{J,N}$:  The curvature
of $D^N$ evaluated at $J$ and on two variations of $J$.
\item $\chi_{A,B}(x) = \tr(A_xJ_xB_x)$: the symbol of the curvature.
\item $\calP_t:\calE_0\to\calE_t$:  Parallel transport along a curve
$J(t)$ in $\calJ$.
\item $Z\subset L^*$: The unit circle bundle in $L^*$.
\item $L^2(Z)_k\cong L^2(M,L^{\otimes k})$: The Fourier decomposition of
$L^2(Z)$.
\item $\calH_J := \bigoplus_{N=1}^\infty\calE_{J,N}$.
\item $\Pi_J: L^2(Z)\to\calH_J$: Orthogonal projector.
\item $\Sigma = \{\, (p, r\alpha_p)\;;\;p\in Z, r>0 \,\}$.
\item $\calS_\sigma = (T_\sigma\Sigma)^o/T_\sigma\Sigma \cong 
T_{\pi(\sigma)} M^{-}$.
\item $e_{J,\sigma}\in\calW_\sigma$: the ground state.
\item $\Upsilon_J(A,B):=
\bigoplus_{N=1}^\infty\Upsilon^{(N)}_J(A,B)$.
\item $V$ a symplectic vector space, typically $V = T_xM$.
\end{itemize}

\end{appendix}

\end{document}